

\magnification=1200

\def\proof{\noindent{\bf Proof. }}
\def\remark{\noindent{\bf Remark }}
\def\example{\noindent{\bf Example }}

\def\O{{\cal O}}

\def\E{{\cal E}}
\def\F{{\cal F}}
\def\G{{\cal G}}
\def\H{{\cal H}}

\def\L{{\cal L}}

\def\P{{\bf P}}
\def\Q{{\bf Q}}
\def\P{{\bf P}}
\def\Z{{\bf Z}}

\def\f{\varphi}
\def\ra{\rightarrow}
\def\raa{\longrightarrow}

\def\har#1{\smash{\mathop{\hbox to .8 cm{\rightarrowfill}}
\limits^{\scriptstyle#1}_{}}}

\def\iso{\cong}

\def\PE{{\bf P}({\cal E})}
\def\PF{{\bf P}({\cal F})}

\def\Ext{{\cal E}xt}
\def\s{{\sl B\v anic\v a sheaf }}
\def\ss{{\sl B\v anic\v a sheaves }}

\def\Diag{\def\normalbaselines{\baselineskip20pt
\lineskip3pt\lineskiplimit3pt}
        \matrix}

{\rightline{December 1993}}

\vskip 1 true cm

{\centerline{{\bf On B\v anic\v a  sheaves and Fano manifolds}}}
\vskip .5 true cm

\centerline{Dedicated to the memory of Constantin B\v anic\v a}
\vskip 1 true cm
\centerline{ Edoardo Ballico and Jaros\l aw A.~Wi\'sniewski}
\vskip 1 true cm

\beginsection Introduction.

Reflexive sheaves are nowadays a common tool to study
projective varieties.
In the present paper we apply reflexive sheaves to study projective morphisms.
Given a projective map $\f: X\ra Y$ and an ample line bundle $\L$ on
$X$ one may consider an associated coherent sheaf $\F:=\f_*\L$ on $Y$.
The knowledge of the sheaf $\F$ allows sometimes to understand
some properties of the variety $X$ and of the map $\f$. This is
a typical way to study cyclic coverings (or, more generally,
finite maps) and projective bundles. In the latter case one may choose
the bundle $\L$ to be a relative $\O(1)$-sheaf so that $X=\P(\F)$.
A similar approach can be applied to study equidimensional
quadric bundles: again, choosing $\L$ as the relative $\O(1)$,
one produces a projective bundle $\P(\F)$ in which $X$ embedds
as a divisor of a relative degree two. Note that, in all the above
examples, if $X$ and $Y$ are smooth then the map $\f$ is flat
and the resulting sheaf $\F$ is locally free. In the present paper
we want to extend the method also to non-flat maps. In particular,
we will consider varieties which arise as projectivizations of
coherent sheaves.
\par

Our motivation for this study was originally two-fold: firstly
we wanted to understand the class of varieties called by
Sommese (smooth) scrolls --- they occur naturally in his adjunction
theory --- and secondly we wanted to complete
a classification of Fano manifolds
of index $r$, dimension $2r$ and $b_2\geq 2$ --- the task which
was undertaken by the second named author of the present paper.
As our understanding of the subject developped we have realised
that many other points and applications
of the theory of projective fibrations
are also very interesting and deserve proper attention.
However, for the sake of clarity of the paper we refrained from dealing with
most of the possible extensions of the theory.
Therefore, in the present paper we will deal mostly with coherent
sheaves whose projectivizations are smooth varieties. This class
of sheaves is related to the class of {\sl smooth sheaves} which
were studied by Constantin B\v anic\v a in one of his late papers.
Thus we decided to name the class of the sheaves studied in the present
paper after B\v anic\v a to commemorate his name.
\par

The paper is organised as follows: in the first two sections we
introduce some pertinent definitions and constructions and subsequently
we examine their basic properties. In particular we prove that \ss
of rank $\geq n$ (where $n$ is the dimension of the base) are locally free,
and subsequently we discuss a version
of a cojecture of Beltrametti and Sommese on smooth scrolls.
In Section 3 we gathered a number of examples which illustrate
some aspects of the theory. From section 4 on we deal with \ss
of rank $n-1$: first we discuss when they can be extended to locally
free sheaves and examine numerical properties of extensions.
In the remaining two sections we apply this to study ampleness
of the divisor adjoint to a \s and then to classify Fano manifolds
of large index which are projectivizations of non-locally free sheaves.
\medskip

\noindent{{\bf Acknowledgements.}}
We would like to thank
SFB 170 {\sl Geometrie und Analysis} in G\"ottin\-gen and
Max-Planck-Institut f\"ur Mathematik in Bonn; parts
of the present paper were prepared at the time we visited
these institutions.
The first named author would like to acknowledge the support from
Italian MURST and GNSAGA while the second named author would
like to acknowledge the support from Polish grant KBN GR54.
\medskip

\noindent{{\bf Notation and assumptions.}} We adopt standard notation,
see Hartshorne's texbook [H1]. We will frequently identify divisors
and line bundles on smooth varieties.
We assume that all varieties are defined over complex
numbers, though the definitions and some results are also valid
for varieties over an algebraically closed field.

\beginsection 1. Projectivization.

First, let us recall the definition
of a projectivization of a coherent sheaf
$\E$ over a scheme $V$, see [G] and [H1] for details.

\noindent {\sl (1.0). }
We start with a local description.
Let $A$ be a noetherian ring and $M$ a finitely
generated $A$-module. We will also usually assume that
the ring $A$ is an integrally closed domain, though
it is not needed for the definitions.
Let $B$ denote the symmetric algebra of $M$
$$B:=Sym(M)=\bigoplus_{m\geq 0}S^m(M)$$
where $S^mM$ is the $m$-th symmetric product
of the module $M$.
The $A$-algebra $B$ has a natural gradation $B_m=S^m(M)$
and we define
${\P}_A(M):= Proj(B)$.
Such defined projective scheme is a generalisation
of the projective space over $A$.
The scheme $\P(M)$ has a natural
affine covering defined by elements
of $M$: for a non-zero $f\in M$ consider
$${\cal D}_+(f)=\{q\in Proj(B): f\notin q\},$$
the scheme ${\cal D}_+(f)$ is then isomorphic to
an affine scheme $Spec(B_{(f)})$, where $B_{(f)}$ denotes
the zero-graded part of the localisation $B_f$
of $B$ with respect to
the element $f$.
The embedding $A=S^0M\subset Sym(M)$
yields a projection map
$$p: {\P}(M)\rightarrow SpecA.$$

\par

Graded modules over $B$ give rise to coherent sheaves
over $\P(M)$. In particular,
on $\P(M)$ there are invertible
sheaves $\O(k)$ associated to graded
$B$-modules $B(k)$, with $B(k)_m=B_{k+m}=S^{m+k}(M)$,
where the sub-index denotes the gradation
shifted by $k$ with respect to the gradation of $B$.
Note that sections of the sheaf $\O(1)$ are isomorphic to
the module $M$.

\medskip

The above local definition of $\P$ allows us to define
projectivization for any coherent sheaf $\E$:
If $Sym\E:= \bigoplus_{m\geq 0} S^m\E$ is the symmetric algebra
of sections of coherent sheaf $\E$ over a normal
variety $V$ then we define
$$\P(\E):= Proj_V(Sym\E).$$
The inclusion $\O_V\iso S^0\E\ra Sym\E$ yields the projection
morphism $p: \P(\E)\ra V$. We will always assume that the
morphism $p$ is surjective, or equivalently, that the support of $\E$
coincides with $V$. The local definition of $\O(1)$
gives rise to a globally defined invertible sheaf and thus
over $\P(\E)$ there exisits
an invertible sheaf $\O_{\P(\E)}(1)$ such that $p_*\O(1)=\E$.

\medskip
In the present section we want to understand some basic properties
of this construction. The first one is about
irreducibility.

\proclaim Lemma 1.1. If $\P(\E)$ is an integral scheme
then $\E$ is torsion-free.

\proof The assertion is local. Note that
$\O(1)$ is locally free of rank
1 on an integral scheme and therefore it has no torsions.
Consequently $\E$,
being locally the space of sections of $\O(1)$,
is torsion-free.
\medskip

The converse of the above lemma is not true, see the example (3.2).

\medskip
Therefore, from now on we will assume that
all the sheaves whose projectivizations we will consider
are torsion free.
\medskip

\proclaim Lemma 1.2. {\rm (cf. [H2, 1.7])}
Let $\E$ be a torsion-free sheaf over
a normal variety $Y$ and let $p: \P(\E)\ra Y$
be the projectivization
of $\E$. Assume that $\P(\E)$ is a normal variety and
no Weil divisor in $\P(\E)$ is contracted to
a subvariety of $Y$ of codimension $\geq 2$. Then
the sheaf $\E$ is reflexive.

\proof
We claim that the sheaf $\E$ is normal (in the sense of [OSS, II,1]
or [H2]). This is because
any section of $\E$ over open subset $U\setminus D$ of $Y$,
where $D$ is of codimension $\geq 2$,
is associated to a section of $\O(1)$
over $p^{-1}(U\setminus D)$. This, however, extends
uniquely over $p^{-1}(U)$ because $\P(\E)$ is normal
and $p^{-1}(D)$ is of codimension $\geq 2$ (c.f.~[H2, 1.6]).

\medskip

The above argument works for any projective surjection $\f: X\ra Y$ of normal
varieties. If $\f$ contracts no Weil divisor on $X$ to a codimension
$\geq 2$ subset of $Y$ then a push-forward $\f_*\F$ of any reflexive
sheaf $\F$ on $X$ is reflexive on $Y$, see [H2, 1.7].
This is used in the following

\proclaim Lemma 1.3. Let $\E$ be a reflexive
sheaf over a normal variety $Y$ satisfying
the assumptions of the previous lemma. Then
$$p_*(\H om(\Omega_{\P(\E)/Y},\O(-1))) \iso \H om(\E,\O_Y).$$

\proof Note that the isomorphism is true
if $\E$ is locally free (one can use relative Euler sequence
to prove it).
The sheaf $\H om(\Omega_{\P(\E)/Y},\O(-1))$ is reflexive
as a dual on a normal variety. Then, similarly as above
we prove that its push-forward is reflexive as well.
Thus we have isomorphism of the two reflexive sheaves
defined outside of a codimension 2 subset of $Y$.
Therefore the sheaves are isomorphic.
\medskip

We will need the following

\proclaim Lemma 1.4.
Let $(A,m)$ be a regular local ring which is an
algebra over its residue field $k=A/m$. Assume that $M$ is an
$A$-module which is not free and which comes from an exact sequence
$$0\raa A\mathop{\raa}^s A^{r+1}\raa M\raa 0.$$
Let us write $s(1)=(s_0,\dots, s_r)$ where $s_i\in m\subset A$.
Then $\P_A(M)$ is regular if and only if the classes of elements
$s_0,\dots,s_r$ are $k$-linearly independent in $m/m^2$.

\proof The ideal of $\P(M)$ in
$\P_A^{r}=Proj(A[t_0,\dots, t_r])$
is generated by an element $\sum s_i t_i$. Therefore,
in an affine subset $U_0=Spec A[t'_1,\dots, t'_r]$
(where $t'_i=t_i/t_0$)
its equation is $s_0+\sum s_i t'_i = 0$. Thus, $\P(M)$ is smooth
at $t'_1=\dots=t'_r=0$ if and only if $s_0$ is non-zero in $m/m^2$.
The above argument can be repeated for any $k$-linear transformation
of coordinates in $m/m^2$ which proves
that $s_0,\dots,s_r$ are linearly
independent in $m/m^2$.
\medskip

\beginsection 2. B\v anic\v a sheaves, first properties.

In one of his last papers [B], Constantin B\v anic\v a considered
a special class of reflexive sheaves.

\proclaim Definition 2.0. A reflexive sheaf $\E$ of rank $r$
over a smooth variety $V$ is called smooth if $\Ext^q(\E,\O)=0$
for $q\geq 2$ and $\Ext^1(\E,\O)=\O_v/(t_1,\dots,t_{r+1})$
for some choice $(t_1,\dots,t_n)$ of regular parameter system
at a point $v$ of singularity of $\E$.

Smooth sheaves are convenient for studying subvarieties
of smooth varieties, see also [H2], [BC] and [HH].

In the present paper we will deal with another
special class of coherent
sheaves over normal varieties.
As it will be seen this class is a generalisation
of the one studied by B\v anic\v a and therefore
we name these sheaves after him.

\proclaim Definition 2.1. A coherent
sheaf $\E$ of rank $r\geq 2$ over a normal variety
$Y$ is called \s if its projectivization is a smooth variety.

The assumption on smoothness of the projectivization is very strong
as the following lemma shows.

\proclaim Lemma 2.2. If $\E$ is a \s then it is reflexive and moreover
the map $p: \P(\E)\ra Y$ is an elementary, or extremal ray contraction.
Furthermore $\O_{\P(\E)}(1)$
is $p$-ample and generates $PicX$ over $PicY$
so that we have a sequence
$$0\raa PicY\raa Pic\P(\E)\raa {\bf Z}[\O(1)]\raa 0.$$
Moreover every Weil divisor on $Y$ is Cartier and,
in particular, $Y$ is Gorenstein.

\proof
First, note that since $\P(\E)$ is irreducible, $\E$ is torsion-free (1.1).
To prove that $p$ is an elementary contraction note that
every fiber of $p$ over a point $y\in Y$
is a projective space
$\P(\E_y\otimes k(y))$ (where $k(Y)$ denotes the residue field).
Taking a line in a generic fiber and deforming it,
we obtain a non-trivial curve
in a special fiber, too (actually a line),
therefore all curves contracted are numerically
proportional, hence
$p$ is an extremal ray contraction in the sense of
Mori theory and consequently:
$\O(1)$ is $p$-ample.
Moreover there is an exact sequence
$$0\raa PicV\raa Pic\P(\E)\raa {\bf Z}[\O(1)]\raa 0.$$
The sequence is exact even at the last place
because $\O(1)$ has intersection 1 with
a line in the fiber.
If a prime Weil divisor in $\P(\E)$ is
contracted to a proper subset of $Y$
then it has trivial intersection with curves contracted by $p$
and thus it is a pull-back of a Cartier divisor from $Y$.
Now the reflexivity of $\E$ follows because of
lemma (1.2). The last assertion of the lemma follows similarly:
the inverse image of a Weil divisor from $Y$ is Cartier on $X$
and has intersection 0 with curves contracted by $p$ and thus
it is a pull-back of a Cartier divisor from $Y$.
\medskip

\noindent {\sl (2.3).}
One motivation to study B\v anic\v a sheaves comes from {\sl smooth scrolls}
which are defined by Sommese as follows:
A pair $(X,\L)$ consisting of a smooth variety $X$ and an
ample line bundle $\L$ is called a {\sl scroll} if
there exists a
morphism $p: X\ra Y$ onto a normal variety $Y$ of smaller dimension
such that
$K_X\otimes\L^{\otimes(dimX-dimY+1)}$
is a pull-back of an ample line bundle from $Y$.
\par

A smooth scroll is over a general point a projective bundle,
this follows from Kodaira vanishing and Kobayashi---Ochiai
characterisation of the projective space.
Obviously, projective bundles and, more generally, projectizations
of coherent sheaves are examples of smooth scrolls.
Conversly, if all fibers of the map $p$ are of the same dimension
then the scroll is a projective bundle, [F1, 2.12] and [I].
We have also examples of scrolls which do not belong to any
of these two classes; their fibers may be Grassman varieties
of large dimension with respect to the dimension of
a general fiber, see the example (3.2).
If we assume that the smooth scroll is a projectivization of a sheaf,
the dimension of special fibers can not jump so much:

\proclaim Lemma 2.4.
Let $\PE\ra Y$ be a projectivization
of a rank-$r$ \s. Let $F$ be a fiber of dimension
$>r-1$. Then $dimF\leq dimY$.

\proof Let $\Pi\iso\P^{r-1}\subset F\iso \P^k$ be a specialization
of a general fiber.
We have then a sequence of normal bundles
$$0\raa N_{\Pi/F}\iso \O(1)^{k-r+1}\raa N_{\Pi/\PE}\raa
N_{(F/\PE)\vert\Pi}\raa 0.$$
Since $N_{\Pi/\PE}$ is a specialization of a trivial bundle
it has a trivial total Chern class, therefore
$$c_t(N_{(F/\PE)\vert\Pi})=\bigl(c_t(\O(1))\bigr)^{r-k-1}.$$
Consequently,
$$n+r-1-k=rank(N_{F/\PE})\geq dim\Pi=r-1$$
and the inequality follows.
\medskip
On the other hand, if we assume that the jump of the dimension
of fibers in a scroll is small then we can apply Theorem 4.1 from [AW]
to get the following

\proclaim Proposition 2.5. Let $(X,\L)$ be a smooth scroll. Assume that
for any fiber $F$ of the map $p: X\ra Y$ it holds $dim F\leq dimX-dimY$.
Then $Y$ is smooth and $X=\P(p_*\L)$  (so that $p_*\L$ is
a B\v anic\v a sheaf.
Moreover, if $dimX\geq 2dim Y$ then $p$
is a projective bundle.

For smooth scrolls which are projectivization
of sheaves there holds a conjecture of Beltrametti and Sommese;
namely we have

\proclaim Theorem 2.6.
Let $\E$ be a \s of rank $r$ over a normal variety $Y$.
If $r\geq dimY$ then $Y$ is smooth and $\E$ is locally free.
If $r=dimY-1$ then $Y$ is smooth.

\proof The first part follows immediately from Lemma 2.4
and Fujita's result [F, 2.12].
Then, the second part follows then from 2.5.
\medskip

Using the above Proposition 2.5 and Remark 4.12 from [AW] we get
the following

\proclaim Lemma 2.7.
Let $\E$ be a \s of rank $r$ over a normal variety $Y$.
If $r\geq dimY-1$ or, if for any point $y\in Y$,
 $dim_k \E_y\otimes k(y)\leq r+1$,
then $Y$ is smooth and locally $\E$ is a quotient of a trivial
sheaf by a rank-1 subsheaf, that is, we have a sequence
$$0\raa \O_{Y,y}\raa \O_{Y,y}^{r+1}\raa \E_y\raa 0.$$

If we now combine lemmata (1.4) and (2.7) we get the following

\proclaim Corollary 2.8. Any smooth sheaf
(in the sense of B\v anic\v a)
is a B\v anic\v a sheaf. If a \s $\E$ over a normal variety $Y$
satifies the condition
$$\hbox{ for any }y\in Y :
dim_k \E_y\otimes k(y)\leq rank\E+1$$
then it is smooth.

\beginsection 3. Examples.

The simplest examples of scrolls are projective bundles.
In particular,
if the base $Y$ is smooth then any locally free sheaf is a \s.
Also, if a locally free sheaf $\F$ over a smooth $Y$
is spanned by global
sections then a general section $s$
of $\F$ will yield a \s as a quotient:
$$0\raa \O\mathop{\raa}^s\F\raa\E\raa 0.$$
The singular set of $\E$ coincides with the zero locus
of the section $s$. The local condition on the sheaf $\E$
to be \s is described in Lemma (1.4).
\par
More generally, we can consider arbitrary morphisms
of vector bundles over smooth base

\proclaim Lemma 3.1. {\rm (cf.~[B, Thm.~2])}
Let $\F$ and $\G$ be locally free sheaves over
a smooth variety $Y$ of rank $f$ and $g$, respectively.
Assume that $f\geq g+2$ and the sheaf $\H om(\G,\F)$ is spanned
by global sections. Then, for a generic
$\sigma\in Hom(\G,\F)$ we have
an exact sequece
$$0\raa \G\mathop{\raa}^{\sigma}\F\raa\E\raa 0$$
with the quotient $\E$ being \s of rank $f-g$.

\proof We have a natural isomorphism (see [H1, II.5])
$$Hom_{\P(\F)}(p^*\G,\O(1))\iso Hom_Y(\G,\F).$$
The zero locus of a section
$$\sigma\in Hom(p^*\G,\O(1))=H^0(\P(\F),p^*\G^*(1))$$
coincides with the projectivization of the cokernel $\E$
of the map $\sigma: \G\ra\F$ embedded into $\P(\F)$
by the map associated to the epimorphism $\F\ra\E$.
Therefore the lemma follows from Bertini theorem.
\medskip

Not all scrolls arise as the projectivizations
of sheaves.
\medskip

\example {\bf 3.2.} Consider the Grassmann variety $G(2,n)$ of linear planes
of a given linear space $W$ of dimension $n$.
Over $G(2,n)$ we have the universal
quotient bundle ${\cal Q}$ whose projectivization is a flag variety
$$F(1,2,n)=\{(x,l)\in \P^{n-1}\times G(2,n): x\in l\}$$
with a projection onto $\P^{n-1}$.
The projection has a natural structure of
$\P^{n-2}$-bundle. Now take a bundle ${\cal Q}\oplus \O$,
its projectivization $q: \P({\cal Q}\oplus \O)\ra G(2,n)$
maps onto $\P^{n}$, $p: \P({\cal Q}\oplus \O)\ra \P^{n}$
 so that all fibers but one
are isomorphic to $\P^{n-2}$. The exceptional fiber (call it $F_0$)
is associated to the $\O$-factor of the bundle $q$ and
it is isomorphic to $G(2,n)$, so that it is of
dimension $2(n-2)$. One checks easily that the variety has a structure
of a smooth scroll, however it is not a projectivization of a sheaf
as the special fiber is not a projective space (for $n\geq 4$).
Let us consider the sheaf
$\F:= p_*q^*\O_{G(2,n)}(1)$ where $\O_{G(2,n)}(1)$
is the positive generator of $PicG(2,n)$.
The sheaf $\F$ is locally free
outside one point where it has a fiber isomorphic to $\Lambda^2 W$.
It is not hard to check that it is
reflexive though its projectivization
is a reducible variety consisting of two components: the dominant one
which is the original scroll and the special fiber $\P(\Lambda^2 W)$,
the fiber $F_0$ embedded in $\P(\Lambda^2 W)$ via Pl\"ucker embedding.

\medskip

If we allow the projectivization have some singularities,
even mild ones, some of the statements from the previous section
are not true (e.g.~2.4).

\example {\bf 3.3.} (Sommese [S, 3.3.3])
Take a smooth surface $S$ and blow it up
$\beta: S'\ra S$ at a point $s\in S$.
Let $E$ denote the exceptional divisor.
Let $L$ be a pull-back to $S'$
of an ample line bundle from $S$. We may assume (possibly
replacing $L$ by its power) that
$L\otimes\O(-E)$ is ample and spanned on $S'$.
Over $S'$ we consider
a projective bundle $p':\P\bigl(L\oplus (L\otimes\O(-E))\bigr)\ra S'$.
The $\O(1)$-sheaf on the projectivization
is clearly nef and ample outside the inverse image of
$E$. The unique curve with which the $\O(1)$-sheaf has trivial intersection
is the section of the projective bundle over $E$
(a smooth rational curve)
associated to the splitting $\O\oplus\O(1)\ra\O$.
The smooth rational curve is easily seen to have normal bundle
$\O_{\P^1}(-1)\oplus\O_{\P^1}(-1)$ and it can be contracted
to an isolated singular point by the morphism
coming from the evaluation of $\O(1)$
(since the bundle $L\oplus (L\otimes\O(-E))$ is spanned).
 The singularity is Gorenstein
since the canonical bundle on the projectivization
has intersection 0 with the contracted curve.
By $V$ let us call the resulting 3-fold
obtained by contracting the section to a point.
The 3-fold $V$ maps onto $S$ and the map makes $V$ a scroll.
There exists a unique exceptional fiber of the scroll which is
isomorphic to $\P^2$ and which contains the singular point.
On the other hand, the threefold $V$ can be described
as a projectivization of a sheaf $\E:=(\beta\circ p')_*\O(1)$,
and it is not hard to see that the singularity of $\E$
at the point $s$ is of the type $\O\oplus J_s$, where
$J_s$ is the ideal of the point $s$.
\medskip

In the present paper we will also deal with Fano manifolds arising
as projectivization of sheaves. We have:

\proclaim Lemma 3.4. Let $\E$ be a \s over a normal variety
$Y$. Assume that a singular set of
$\E$ is of dimension $\leq 1$ or $\rho (Y)=1$.
If $\PE$ is a Fano manifold then $-K_Y$ is ample,
that is, $Y$ is a Gorenstein Fano variety.

\proof The argument is similar as in the proof of [W, 4.3], compare also
with [SW, 1.6] and [KMM].
We are only to prove that
$$p^*(-K_Y)=\bigl(-K_{\P(\E)}\bigr)+ \O_{\P(\E)}(-rank\E)+ p^*(-det\E)$$
has positive intersection with
any extremal rational
curve $C$ in $\PE$ not contracted by $p$.
We claim that the curve $C$ may be chosen so that
$p(C)$ is not contained in the singular locus
of $\E$. Indeed, if it were, then the whole
locus of the ray ${\bf R}^+[C]$ would be contracted by $p$
to a set of dimension 1, thus all fibers
of the contraction of the ray would be of dimension
1, hence the locus would be a divisor, [W1, 1.1]. This, however,
contradicts the fact that $p$ contracts no Weil
divisors to set of codimension $\geq 2$.
Once the curve $C$ is assumed not to be in the
singular locus of the map $p$ we conclude as in
[SW, 1.6], or as in [KMM].

\medskip

\example {\bf 3.5.} The assumption on the singular set of $\E$ is
indespencible. Let
$$Y:=\P(\O_{\P^2}(1)^{\oplus 3}\oplus \O_{\P^2}).$$
Then $-K_Y=4\eta$ where $\eta$ denotes the relative
$\O(1)$ of the projectivization over $\P^2$. The line bundle $\eta$
is spanned but not ample, so $Y$ is not Fano.
The morphism associated to $\vert \eta\vert$
contracts to a point the unique section of $Y\ra \P^2$
associated to the $\O$-factor, call this set $Z$. Let
$H$ be the pullback of the hyperplane from $\P^2$ to $Y$.
The line bundle asociated to $\eta-H$ is spanned
off $Z$ by three sections. Thus we have a morphism
$\O_Y^{\oplus 3}\ra\O_Y(\eta-H$) which yields a sequence
$$ 0\raa\O_Y(-\eta+H)\raa\O_Y^{\oplus 3}\raa\E\raa 0$$
with a rank-2 sheaf $\E$ which is free outside
$Z$. The variety $\PE$ is the coincidence variety
of divisors from the linear system $\vert \eta-H\vert$
and each one of the divisors is
isomorphic to $\P(\O^{\oplus 2}_{\P^2}\oplus\O_{\P^2}(-1))$. Therefore
$\PE$ is smooth. Moreover
$$-K_{\PE}=\O_{\PE}(2)\otimes p^*(\O_Y(3\eta+H)).$$
Therefore $\PE$ is a smooth Fano variety.

\beginsection 4. Extensions to locally free sheaves, nefness

We want to find conditions to realise globally the projectivization
of a \s as a divisor in a projective bundle. For simplicity
we introduce the following definition.

\proclaim Definition 4.0. We say that a coherent sheaf $\E$
over a normal variety $Y$
extends to a locally free sheaf $\F$ if there exists a sequence
of $\O_Y$-modules
$$0\raa \O\mathop{\raa}^s\F\raa\E\raa 0.$$

In other words, $\E$ is obtained by dividing $\F$
by a non-zero section $s$. The singular locus
of $\E$ coincides with the zero locus of $s$.
Alternatively, $\PE$ is a divisor in $\P(\F)$
from the linear system $\vert\O_{\P(\F)}(1)\vert$.
\medskip

In the present section we will also discuss numerical
properties of coherent sheaves.
Let us recall that a sheaf $\E$ is ample (resp.~nef)
if $\O(1)$ is ample (resp.~nef) on $\PE$;
this makes sense also if we multiply $\E$ by a $\Q$-divisor.
\medskip

For a coherent sheaf $\E$ by $\E^*$ we will denote its dual
sheaf $\H om (\E,\O)$.

\proclaim Lemma 4.1.
Let $\E$ be a B\v anic\v a sheaf of rank $n-1$ over a smooth projective
variety $Y$ of dimension $n$.
If $H^2(Y,\E^*)=0$ then $\E$ extends to a localy free sheaf; in
particular $\E\otimes\L^{-m}$ extends for $\L$ an ample line bundle
and $m\gg 0$.
\par

\proof Because of (2.7) we know that
the extension exists locally. To prove the existence
of a global extension consider the spectral sequence relating
local $\E xt$ and global $Ext$.
Then we have the following exact sequence
$$H^1(Y,\H om(\E,\O))\ra Ext^1_Y(\E,\O)\ra
H^0(Y,\E xt^1(\E,\O))\ra H^2(Y, \H om(\E,\O)).\leqno (4.1.1)$$
The support of $\Ext^1(\E,\O)$ consists of isolated points
of singularity of $\E$.   For any such point $y$,
$\Ext^1(\E,\O)_y\iso\O_y$ and the unit represents the extension
to a free module. Thus the vanishing of $H^2(Y,\E^*)$ yields the existence
of an extension in $Ext^1(\E,\O)$ to a locally free sheaf.

Therefore, frequently we will be interested
in the vanishing of the latter term
in the sequence (4.1.1). To this end  we have.

\proclaim Lemma 4.2.
Let $\E$ be a B\v anic\v a sheaf of rank $n-1$ over a smooth projective
variety $Y$ of dimension $n$.
Assume that $\L$ is an ample line bundle on $Y$ and let $H\in
\mid\L\mid$ be a smooth divisor which does not meet the singular
set of $\E$.
If, for $k\geq 1$ and $i=1,\ 2$
$H^i(H,(\E^*\otimes\L^k)_{\vert H})=0$
then $\E$ extends.

\proof
The vanishing of $H^2(Y,\E^*)$ follows from the vanishing of cohomology
on $H$ which, beause of the divisorial sequence for $H$, implies that
$$H^2(Y,\E^*\otimes\L^k)=H^2(Y,\E^*\otimes\L^{k+1})$$
for $k\geq 0$.
\medskip

On the other hand, the non-vanishing of $H^2(Y,\E^*\otimes\L^k)$
for $k\ll 0$ can be used to estimate $c_n(\E)$, that is, the number of
singular points of $\E$. The following lemma was suggested to us
by Adrian Langer whom we owe our thanks for finding a mistake in a previous
version of this paper.

\proclaim Lemma 4.2.1.
Let $\E$ be a B\v anic\v a sheaf of rank $n-1$ over a smooth projective
variety $Y$ of dimension $n$ and let $\L$ be an ample line bundle
over $Y$. Then, for $k\gg 0$, we have
$H^2(Y,\E^*\otimes\L^{-k})=c_n(\E)$.
\par

\proof
We have a global duality [H1, III.7.6]:
$$Ext_Y^i(\E\otimes\L^k,\O_Y)\iso H^{n-i}(Y,\E\otimes\L^k\otimes K_Y)^*$$
and the latter term vanishes for $k\gg 0$ and $i<n$.
Therefore the spectral sequence relating $Ext$ and $\E xt$ converges
to a trivial one. This yields that $H^0(Y,\E xt^1(\E\otimes\L^k,\O_Y))=
H^2(Y,\E^*\otimes\L^{-k})$ (c.f.~4.1.1) and we are done.
\medskip

Making similar argument as is the proof of 4.2 we get the following
\proclaim Corollary 4.2.2.
Let $\E$ be a B\v anic\v a sheaf of rank $n-1$ over a smooth projective
variety $Y$ of dimension $n$.
Assume that $\L$ is an ample line bundle on $Y$ and let $H\in
\mid\L\mid$ be a smooth divisor which does not meet the singular
set of $\E$.
If, for any $k\in\Z$ and $i=1,\ 2$ the groups
$H^i(H,(\E^*\otimes\L^k)_{\vert H})$ vanish
then $\E$ is locally free.

\medskip

We will need also the following version of the lemmata 4.1 and 4.2 for
arbitrary sheaves with isolated singularities.

\proclaim Lemma 4.3.
Let $\E$ be coherent sheaf with isolated singularities
over a smooth variety $Y$, $dimY\geq 3$.
Let $\L$ be  an ample line bundle over $Y$ and let $H\in
\mid\L\mid$ be a smooth divisor which does not meet the singular
points of $\E$. Then:
\item{(a)} if $\E$ extends to a locally free sheaf
then $\E\otimes\L^{-1}$ extends as well,
\item{(b)} if $\E\otimes\L^{-1}$ extends to a locally free sheaf
and $H^2(Y,\E^*)=0$
then also $\E$ extends to a locally free sheaf.
\item{}
\par

\proof
Consider a divisorial sequence associated to $H\in \vert\L\vert$:
$$0\raa\O\raa\L\raa\L_H\raa 0.$$
The morphism $\O\ra\L$ from this sequence yields
a commutative diagram with exact rows and columns coming from multiplying
by a section defining $H$:
$$\Diag{
Ext^1_Y(\E,\O)&\ra&H^0(Y,\E xt^1(\E,\O))&\ra&H^2(Y,\E^*)\cr
\downarrow&&\downarrow&&\downarrow\cr
Ext^1_Y(\E,\L)&\ra&H^0(Y,\E xt^1(\E,\L))&\ra&H^2(Y,\E^*\otimes\L)}$$
On the other hand we know that
$$\E xt^i(\E,\L)\iso\E xt^i(\E,\O)\otimes\L$$
so that, because the singularities of $\E$ are isolated
and $H$ does not meet them, the vertical map
in the center is an isomorphism.
Therefore an extension in $Ext^1(\E,\O)$ which gives a locally
free sheaf will be mapped by the left-hand-side vertical map
to an extension in $Ext^1(\E,\L)$ which produces a locally free sheaf, too.
This proves (i). To get (ii) we make a similar agument, but this
time applying vanishing of $H^2(Y,\E^*)$ to lift a local exension
to a global one.

\medskip

Now we want to compare ampleness
and nefness of a rank $r$ \s $\E$
with the same properties of a locally
free sheaf $\F$ in whose projectivization $\E$ is embedded.
Therefore, let us assume that $\E$ extends to $\F$,
that is, we have the sequence 4.0.
Obviously, if $\E$ is nef then also $\F$ is nef. As for the
ampleness we have the following.

\proclaim Lemma 4.4.
Let $\E$ and $\F$ be coherent sheaves on a smooth variety $Y$
satisfying the above assumptions.
Assume moreover that $c_1Y-c_1\F$ is nef and that $\E$ is ample.
Then $\O_{\P(\F)}(1)$ is semiample, that is $\vert\O_{\P(\F)}(m)\vert$
is base point free for $m\gg 0$. The exceptional
set $E$  of the morphism given by
$\vert\O_{\P(\F)}(m)\vert$, $m\gg 0$, if non-empty,
contains all sections $Y\supset G\ra\P(\F_{G})$
associated to a splitting
$$\F_{\vert G}\ra \O_G\ra 0$$
of the sequence (4.0) over any
closed $G\subset Y$ of positive dimension.
Moreover $p$ maps $E$ finite-to-one into
$Y\setminus sing\E$.

\proof The line bundle $\O_{\P(\F)}(1)$ is nef and big.
Since $K_{\PF}=\O(-r-1)\otimes p^*(K_Y+det\F)$,
it follows that $\O_{\PF}(m)\otimes K_{\PF}^{-1}$ is nef and big
for $m\gg 0$. Therefore, by the Kawamata-Shokurov contraction theorem
$\O_{\PF}(m)$ is semiample. The morphism defined by
$\vert\O_{\P(\F)}(m)\vert$ is birational and its
exceptional set does not meet $\PE\subset \PF$
(because the divisor $\PE\subset\PF$ has positive intersection
with any curve meeting it).
If the sequence (4.0) splits over a positive-dimensional
set $G\subset Y$ then, clearly, the unique section
of $\F$ over $G$ is contained in $E$.
And clearly $p(E)\cap sing(\E)=p(E\cap\PE)=\emptyset$.

\proclaim Corollary 4.5. In the above situation,
if $G\subset Y$ is not contained in $p(E)$ then
$$Ext^1_G(\E_G,\O_G)\ne 0.$$

We will also need the following.

\proclaim Lemma 4.6. Let $\E$ be an ample reflexive
sheaf on a normal variety $Y$.
Assume that some twist of $\E$ extends
to a locally free sheaf so that
we have a sequence
$$0\raa\L\raa\F\raa\E\raa 0$$
with $\F$ locally free and $\L$ a line bundle.
Let $C\subset Y$ be a rational curve
which is not contained in the singular locus of $\E$. Then
$$C.det\E\geq rank\E + \hbox{ number of singular points of }\E
\hbox{ on }C.$$

\proof First, let us note that, in the above situation,
$$C.det\E = (\O_{\PE}(1))^r.p^{-1}(C)$$
where $r$ is the rank of $\E$ and $p:\PE\ra Y$ is the projection.
Indeed, the formula is correct for projective bundles
and is preserved for divisors in them which meet the cycle $p^{-1}(C)$
at the expected dimension.
The cycle $p^{-1}(C)$ consists of ``vertical'' components
over singular points (each being a projective space)
and of the dominant component over $C$ which is a projective
bundle with a fibre $\P^{r-1}$. From the classification
of bundles over $\P^1$ it follows that the latter component
brings to the intersection at least $r$ and therefore
the inequality follows.

\beginsection 5. Adjunction.

In the present section we compare the determinant, or the first
Chern class of a \s with the canonical sheaf of the variety
over which the sheaf is defined. In case of locally free sheaves
the question was considered in [YZ], [F2] and [ABW].

\proclaim Theorem 5.1.
Let $\E$ be a \s of rank $r$ over a smooth
variety $Y$ of dimension $n=r+1\geq 3$.
Assume that $\E$ is ample and moreover that it is
not locally free. Then
\item{(1)} $K_Y+c_1\E$ is nef unless
$Y\iso \P^n$ and $\E$ is a quotient of a decomposable
sheaf:
$$0\raa\O\raa\O(1)^{\oplus n}\raa\E\raa 0.$$
\item{(2)} if $K_Y+c_1\E$ is nef then it is also
big unless
\itemitem{(2.1)} $Y$ is Fano and $K_Y+c_1\E=0$, or
\itemitem{(2.2)} $Y$ has a structure of a projective
bundle $\pi:Y\ra B$ over a smooth curve $B$ and $\E$ fits
into a sequence
$$0\raa\O\raa\pi^*\G\otimes\O_Y(1)\raa\E\raa 0$$
where $\G$ is a rank-$n$ vector bundle over $B$
and $\O_Y(1)$ a line bundle whose restriction
to any fiber of $\pi$ is $\O(1)$;
\item{(3)} if $K_Y+c_1\E$ is nef and big then it is
also ample unless there exists a birational map $\pi: Y\ra Y'$
supported by $K_Y+c_1\E$ onto a smooth variety $Y'$
which blows-down disjoint exceptional divisors $E_i\iso\P^{n-1}$,
such that $E_i\cap sing\E=\emptyset$.
On $Y'$ there exists an ample  \s $\E'$ such that
$\E\iso\pi^*\E'\otimes\O_Y(-\sum E_i)$ and $K_{Y'}+c_1\E'$
is ample.

\remark The case (2.1) of the theorem will be discussed thoroughly
in the subsequent section. In particular, it will be shown that
$Y$ is either  a projective space or a smooth quadric.
\medskip

\noindent
{\bf Proof of the theorem.} If  $K_Y+c_1\E$ is not nef, then according
to the cone theorem of Mori, there exists an extremal ray of $Y$
which has negative intersection with this divisor.
The length of the ray is at least $n$ so that its locus
coincides with $Y$, see [I, 0.4] or [W1, 1.1].
Therefore there exists a rational curve from
the ray meeting the singular locus of $\E$. Because of (4.6),
these curves have
intersection at least $n+1$ with $-K_Y$. Consequently, by an
argument on deformation of curves passing through a point (see e.g.~[W1]),
$PicY=\Z$ and  we compute easily that $K_Y=(n+1)(K_Y+det\E)$
and therefore by a theorem of Kobayashi-Ochiai $Y\iso\P^n$.
The restriction of $\E$ to a generic
hyperplane $H\subset \P^n$ is an ample vector bundle and $c_1(\E_H)=n+1$,
therefore we see that $\E_H\iso T\P^{n-1}$ or
$\E_H\iso\O(2)\oplus\O(1)^{\oplus(n-1)}$, the latter possiblity ruled out
because of 4.2.2.
In the former case, we use 4.2 to produce an extension of $\E$ to a
locally free sheaf $\F$; the only possible non-trivial extension on $H$
leads to a decomposable bundle $\O_H(1)^{\oplus n}$ so the bundle
$\F$ is decomposable as well (see e.g.~[OSS]).

For the remaining cases the argument is similar. Assume first,
that $K_Y+c_1\E$ is nef but not ample. Therefore there exists
a ray of $Y$ having intersection 0 with the divisor. The length
of the ray is at least $n-1$, [W1, 1.1]. If the contraction of the ray
is birational then it is actually divisorial and the ray
has to have length $n-1$. In this case,
however, the exceptional locus can not meet the singular locus
of $\E$ because then we would find out (4.6)
that the length is actually $n$ which contradicts [I, 0.4].
Consequently, $\E$ is a vector bundle in a neighbourhood of
$E$ and the arguments from [ABW, 2.4] apply to conclude the
description of the blow-down morphism and the sheaf $\E$
as in the case (3) of the theorem.

If the contraction of the ray in question is of fiber
type then a fiber containing a singular point of $\E$
has to be a divisor (again, since $-K_Y.C\geq n$
for any rational curve passing through the singular point).
Thus the contraction is either to a point (which is the case
of (2.1)) or onto a smooth curve $B$. In the latter case
we consider fibers which do not contain singularities
of $\E$ and as in [ABW, 2.2] we prove that the fibers are projective
spaces. Similarly, we conclude that $Y$ has a structure of
a scroll $\P(\G)\ra B$ over the curve and $\E$ restricted to a general
fiber $F$ of the contraction is isomorphic to  $T\P^{n-1}$.
To complete the description of $\E$ we choose a smooth divisor $H\subset X$
which is a hyperplane in each fiber of the scroll and which does not meet
the singular set of $\E$; the restriction of $\E$ to any fiber
of $H\ra B$ is then $T\P^{n-2}\oplus \O(1)$.
After twisting $\E$ by a pull-back of a negative line bundle from
$B$ it will satisfy asssumptions of (4.2) on $H$, so that it will extend
to a vector bundle $\F$. The description of $\F$ follows now
easily, since its restriction to a general fiber has to
be isomorphic to $\O(1)^n$.
\par
To conclude the theorem note that the loci of extremal rays can not
meet (because we would have a curve contracted by both contractions)
and therefore the description of the adjoint morphism is as in (2) of
the theorem.

\beginsection  6. Fano manifolds of middle index.

In the present section we want to complete the classification
of Fano manifolds of index $r$ and
dimension $2r$ with second Betti
number $b_2\geq 2$. Let us recall
that a smooth projective variety
$X$ is called Fano if its anti-canonical
divisor $-K_X$ is ample.
The index of the Fano variety is equal to the largest
integer $r$ for which $-K_X\equiv rH$,
for some ample divisor $H$. Such varieties with projective
and quadric bundle structure were studied in [PSW2]
and [W2], respectively. To complete their classification one has to deal
with these which are non-equidimensional scrolls [W2, Theorem I].
\medskip

\noindent {\sl (6.0).}
Our set-up is as follows:
$X$ is a Fano manifold of index $r$ and dimension $2r\geq 6$,
and it is a projectivization of a non-locally free \s $\E$ over a smooth
variety $Y$ of dimension $r+1$. The projection $X\ra Y$ we will
denote by $p$; we may choose $\E:=p_*(\O(H))$,
so that the line bundle associated to $H$ is  $\O_{\P(\E)}(1)$.
The variety $X$ admits also another non-trivial
map (a contraction) with connected fibers $\f: X\ra Z$
onto a normal variety $Z$. In [W2, Thm.~I] it was proved that
all fibers of $\f$ are of dimension $\leq r$ and thus,
because of [AW, Thm.~4.1], $Z$ is smooth and one of the possibilities
occurs:
\item{{\sl (i)}} $dim Z=r+1$ and $\f: X\ra Z$ is a projectivization of a
non-locally free sheaf;
\item{{\sl (ii)}} $dim Z=2r$ and $\f: X\ra Z$ is a blow-down of a
smooth divisor $E$ in $X$ to a smooth subvariety $T\subset Z$,
$dimT=r-1$;
\item{{\sl (iii)}} $dim Z=r$ and $\f: X\ra Z$ is a quadric bundle;
\item{{\sl (iv)}} $dim Z=r+1$ and $\f: X\ra Z$ is a projective bundle.
\medskip

\noindent{\sl (6.1).}
Fano manifolds with projective bundles were studied in [PSW2];
from the classification obtained in that paper
it follows that the last possibility {\sl (iv)} can not occur.
Quadric bundles were studied in [W2] and it had turned out that
two of the quadric
bundles obtained there have also a structure of projectivization
of non-locally free sheaf:
\item{(a)} a divisor of bidegree $(1,1)$ in the product $\P^r\times\Q^{r+1}$,
\item{(b)} a divisor of bidegree $(1,2)$ in the product $\P^r\times\P^{r+1}$.
\par
{}From now on we assume that we are either
in case {\sl (i)} or {\sl (ii)}, which we will
call fibre and divisorial case, respectively.

\par
\medskip
Our arguments are similar to those from [PSW2]:
we will use ``big fibers'' of the map $\f$, that is fibers
of dimension $r$. We know that they are isomorphic to
$\P^r$ and the restriction of $H$ to each of them is $\O(1)$,
see [AW, 4.1].
First we will deal with the case when $\f$ is divisorial.

\proclaim Lemma 6.2. {\rm (cf.~[PSW2 (7.2)])}
Assume that $\f$ is dvisorial. Then
the restriction of $\O(E)$ to a fiber of $p$
is isomorphic to $\O(1)$.
\par

\noindent{\bf Proof.} Assume the contrary.
Let us take a general fiber $F$ of $p$ such
that $E$ restricted to the fiber is a hypersurface of degree $>1$.
Let us take a line in $F$ which is not contained in $F\cap E$;
choose two points $x_1\ne x_2$ such that $x_1,\ x_2\in F\cap E$.
Let $G_i:=\f^{-1}(\f(x_i))$.
We claim that there exists a curve $C\subset Y$, $p(F)\in C$, such that:
\item{(*)} for a general $c\in C$: $\#(p^{-1}(c)\cap(G_1\cup G_2))\geq 2$

Indeed, note first that $dim(G_i)=r$ and $p$ maps $G_i$ onto a divisor
in $Y$. Therefore, if $\f(x_1)\ne\f(x_2)$ we take a curve in
$p(G_1)\cap p(G_2)\ni p(F)$. If $G_1=G_2$ we consider a curve in
a set $\lbrace y: \#(p^{-1}(y)\cap G_1)\geq 2\rbrace\ni p(F)$
which again is of positive dimension.

Now over a generic $c\in C$ we choose a line $L_c$ in $p^{-1}(c)$
such that $L_c$ is not contained in $E$ and $L_c$ meets $G_1\cup G_2$
at at least two points. This way we can construct a ruled surface
over the normalisation of $C$ which is mapped via $\f$ to a two-dimensional
variety and which contains a curve (or curves) contracted to point (or points)
such that it contradicts the following:

\proclaim Sublemma. Let $\pi: S={\bf P}(\E)\rightarrow C$ be a (geometrically)
ruled surface
(a ${\bf P}^1$-bundle) over a smooth curve $C$. Assume that there exists
a map $\varphi: S\rightarrow {\bf P}^N$ such that the image of $\varphi$
is of dimension 2 and $\varphi$ contracts a curve $C_0\subset S$
to a zero-dimensional set.
Then $C_0$ is a unique section of $\pi$ such that $C_0^2<0$.

\noindent{\bf Proof.}
First, we claim that the curve $C_0$ is irreducible.
Indeed, if $C_1$ and $C_2$
were two irreducible components of $C_0$ then $C_1^2<0$, $C_2^2<0$,
$C_1C_2\geq 0$ and $aC_1-bC_2$ would be equivalent
to a multiple of a fiber of $\pi$
for $a$, $b>0$ and thus $(aC_1-bC_2)^2=0$, a contradiction.
Let $i:B\rightarrow C_0\subset S$ be the normalisation. Consider
$\pi_B: S_B:={\bf P}((\pi\circ i)^*(\E))\rightarrow B$ a ruled surface
over $B$ obtained via base change; it has a section $B_0$ which comes
from the epimorphism
$(\pi\circ i)^*(\E)\rightarrow i^*{\cal O}_{{\bf P}(\E)}(1)$.
The section $B_0$ is mapped birationally to $C_0$ under the induced
map of projective bundles $j: S_B\rightarrow S$ and it is a component
of $B_1=j^{-1}(C_0)$ which is contracted by $\varphi\circ j$.
Since $B_1$ is irreducible, it follows that $B_1=B_0$ and
$$1=deg(B_0=B_1\rightarrow C_0)=deg(S_B\rightarrow S)=deg(C_0\rightarrow C)$$
and therefore we are done.
\par

\noindent{\bf Remark.} Note that this argument works also in case
of Lemma 7.2 from [PSW2] to replace the original ``lift-up'' argument
which is incomplete.

\medskip
We continue with the divisorial case:
As an immediate consequence of the preceeding lemma let us note
that the good supporting divisors of $\f$ and $p$
(i.e. pullbacks of ample divisors from the targets of respective maps)
may be chosen to be $H+E$ and $H-E$, respectively.
\par

Let now $M$ be the intersection of a $r$-dimensional fiber
of $p$ with $E$, it follows that $M\iso \P^{(r-1)}$ and
$H_{\vert M}=E_{\vert M}=\O(1)$. Now since
the map $\f$ maps $M$ onto $T$, by a result of
Lazarsfeld [L] it follows that
$T\iso \P^{(r-1)}$. Since $E+H$ is a pullback of a Cartier divisor
$-K_Z/r$ from $Z$ and $(E+H)_{\vert M}=\O(2)$ it follows that
$-K_Z/r$ restricted to $T$ is either $\O(2)$ or $\O(1)$.
In the latter case, however, using the relation
$${-K_Z}_{\vert T}= -K_T-c_1 N^*_{T/Z}$$
we would find out that $c_1 N^*_{T/Z}=0$.
On the other hand, since $H=-E+(E+H)$ is ample on $E$
it follows that $N^*_{T/Z}\otimes\O(-K_Z/r)$ is ample. Thus, if
$c_1(N^*_{T/Z})=0$ and $\O(-K_Z/r)_{\vert T}=\O(1)$
the bundle $N^*_{T/Z}(1)$
would be isomorphic to $\bigoplus \O(1)^{r+1}$,
a contradiction, since its projectivization would not have
a dominant morphism on $Y$ of dimension $r+1$.
A similar argument done if $(-K_Z/r)_{\vert T}=\O(2)$
leads to the situation when $N^*_{T/Z}(2)$ is ample
with first Chern class $\O(r+2)$ and therefore by splitting type
(see e.g.~[W2, 1.9]) $N^*_{T/Z}(2)\iso T\P^{r-1}\oplus\O(1)^2$.
The projectivization of this latter bundle maps with connected
fibers (because they are hyperplanes in fibers of $p$)
onto $Y$. Therefore, we check that the morphism $p$
restricted to $E$ is given by the evaluation
of the bundle $T\P^{r-1}(-1)\oplus\O^2$ and thus we get the following

\proclaim Lemma 6.3. If $\f$ is divisorial then $Y\iso \P^{r+1}$
and any non-trivial fiber of $\f$ is mapped by $p$
isomorphically onto a hyperplane in $\P^{r+1}$.
Moreover $T\iso\P^{r-1}$ and $N_{T/Z}\iso T\P^{r-1}\oplus \O(1)^2$.
\par

Now we deal with the case when both $p$ and $\f$ are of fiber type

\proclaim Lemma 6.4. {\rm (Comparison Lemma [PSW2, 3.1])}
Assume that $\f$ is of fiber type.
Let
$$r_Y:= min\{ -K_Y.C:\hbox{ where } C \hbox { is rational on }X\}.$$
Then $r_Y\cdot H+p^*(K_Y)$ is a good supporting divisor
for ``the other'' contraction $\f$.
\par
The proof of the above lemma in case $r\geq 4$ is identical
as in [PSW2]; for $r=3$ and $\f$ of fiber type the lemma will
also work because
$\f$ has a fiber of dimension $r$, see Remark (3.4) in
[PSW2].

\remark {\bf 6.5.} Note that in the divisorial case we also have
the comparison lemma since the pull-back of
$\O(1)$ from $Y=\P^{r+1}$ to a fiber of $\f$ is
again $\O(1)$.

\proclaim Corollary 6.6. Assume that $\f$ is either
divisorial or of fiber type.
\item{(a)} Let $F$ be an $r$-dimensional fiber
of $\f$. Then $F$ is $\P^r$ and
$\E_F(-1):=(p^*\E_{\vert F})(-1)$ is nef, and $c_1(\E_F(-1))$
is either 1 or 2.
\item{(b)} If $f\iso \P^{r-1}$ is a general hyperplane in $F$ or ---
for $\f$ of fiber type --- a general fiber of $\f$ then
$\E_f(-1):=(p^*\E_{\vert f})(-1)$ is as described in [PSW1, Thm.~1] or
[PSW2, 0.6].

\proof We already noted that $F=\P^r$. The rest is proved
exactly as (5.2) and (5.3) from [PSW2], the case $c_1=0$
ruled out because $\E$ is not locally free.

\proclaim Lemma 6.7. Assume that $\f$ is of fiber type.
Then both $Y$ and $Z$ are isomorphic to $\P^{r+1}$.
\par

\proof
We use notation from 6.6, i.e.~$F$ is a ``big'' fiber of $\f$ while
$f$ is a general fiber of $\f$, or a general hyperplane in $F$.
Let us consider a composition of maps
$$\P(\E_f)\raa\P(\E)\raa Z$$
where $\P(\E_f)\ra\P(\E)$ is induced by the change of the base $p: f\ra Y$.
We claim that the composition is surjective.
Indeed, if this is not the case then $p^{-1}(p(f_1))\cap f_2=\emptyset$
for a sufficiently general choice of $f_1$ and $f_2$,
so that the intersection of cycles $p^{-1}(p(f_1))\cdot f_2$ is zero.
But note that for a general choice of $f_1$
we have $dim(p^{-1}(p(f_1))\cap F)=r-2$ --- because $p(F)$ is ample on $Y$ ---
and thus $p^{-1}(p(f_1))\cap f_2$ is non-empty
of the expected dimension $r-3\geq 0$ for $f_2\subset F$.
\par

Therefore, for $r\geq 5$, looking up through the list from [PSW1],
we find out that the only possibility when $\P(\E_f)$ admits a surjective
map onto a $r+1$ dimensional variety is
$$\E_f(-1)\iso \O^{r+2}/\O(-1)^2.$$
The map $\P(\E_f)\ra Z$ factors through $\P^{r+1}$
and thus $Z=\P^{r+1}$. The reasoning is clearly symmetric
with respect to the change of $Z$ and $Y$ so the lemma is
proved in this case.

For $r$ equal 3 and 4 we have to eliminate some other possibilities
apart of $\E_f(-1)\iso \O^{r+2}/\O(-1)^2$, that is,
possible sheaves $\E_f(-1)$
which occur in the classification [PSW1] such that $\P(\E_f)$ admits
a morphism onto an $r+1$-dimensional variety.
If $r=4$ the other possibility is a sheaf from the sequence
$$0\raa \O\raa\Omega(2)\oplus\O^2\raa \E_f(-1)\raa 0,$$
see [PSW1].
We claim that in this case $H^1(\P^3,\E_f^*(k))=H^2(\P^3,\E^*_f(k))=0$
for $k\geq 2$ and therefore $\E_F(-1)$ extends to a locally free sheaf,
see 4.2, 4.3.
Indeed, the bundle $\E_f^*(k)$ is isomorphic to either
$T\P^3(k-3)\oplus\O(k-1)$ or to ${\cal N}(k-2)\oplus\O(k-1)^2$,
where ${\cal N}$ is a null-correlation bundle on $\P^3$; thus we check the
vanishing easily.
Now, to conclude this case, note that if $\E_F(-1)$ extends to a locally free
sheaf $\F$ then $\F$ is nef on $\P^4$ and with Chern classes
$(c_1,c_2,c_3)=(2,2,0)$, thus checking it with the list from [ibid]
we arrive to a contradiction.

The case $r=3$ (that is $f=\P^2$)
is dealt with similarly: apart of $\E_f(-1)\iso \O^5/\O(-1)^2$
also decomposable bundles and a bundle with Chern classes $(c_1,c_2)=(2,2)$
admit morphism onto 4-dimensional variety, see the main theorem of [SW].
As above we  check a vanishing to claim that $\E_F(-1)$ extends to a locally
free sheaf $\F$ on $F=\P^3$, $\F$ is nef with Chern class $c_1=2$, thus
globally generated, see [PSW1].
But $\P(\F)$ contains $\P(\E)$ which is mapped onto
a 4-dimensional variety, so itself it has to be mapped onto a 5-dimensional
variety.
Again, by [ibid] the only possibilities for $\F$ are $\O^6/\O(-1)^2$
or $\Omega\P^3(2)\oplus\O$, or ${\cal N}(1)\oplus\O^2$,
where ${\cal N}$ is a null-correlation; we are to exclude the latter
two possibilities.
\par
To this end note that $c_3(\Omega\P^3(2))=0$ and thus $\E_F$ is locally free.
Now we can apply an argument from [PSW2, 5.5]: using the relative
Euler sequence (because $\E$ is locally free at $p(F)$) we compute
the total Chern class of $\Omega X_{\vert F}$:
$$c_t(\Omega X_{\vert F})=c_t(p^*(\Omega Y)_{\vert F})\cdot c_t(\E_F(-1)).$$
On the other hand, because of [AW, 4.9, 4.12] $N^*_{F/X}=T\P^r(-1)$
($N_{F/X}$ denoting the normal bundle),
and we compute
$$c_t(\Omega X_{\vert F})=c_t(\Omega \P^r)\cdot c_t(N^*_{F/X})=1-3h+3h^2-h^3$$
and further
$$c_t(p^*(\Omega Y)_{\vert F})=1-5h+11h^2-13h^3,$$
where $h$ denotes the class of a plane in $\P^3$;
in particular $p(F)\cdot c_3(\Omega Y)$ is not divisible by 5.
In our case, however, the integer $r_Y$ from 6.4 is equal to 5
so either $-K_Y$ is divisible by 5 in $PicY$ and then $Y=\P^4$
or $-K_Y$ generates $PicY$. In the latter case the intersection
of any 1-cycle with any divisor would be divisible by 5 (this follows
e.g. from deformation theory), a contradiction. On the other hand
$c_t(\Omega \P^4)=1-5h+10h^2-10h^3+5h^4$ so comparing it with
the above formula for $c_t(p^*(\Omega Y)_{\vert F})$ we arrive to
a contradiction even if $Y=\P^4$.
This completes the proof of 6.7.

\medskip
To conclude the classification we will deal with the case $Y=\P^{r+1}$
in our set-up 6.0 (i)---(iii).
Let $\E_H$ denote the restriction of $\E$ to a general hyperplane
$H=\P^r$ in $\P^{r+1}$; $\E_H(-1)$ is then a nef vector bundle
(see 6.4, 6.5) with $c_1=2$. Looking up through the list from
[PSW1] we get the following possibilities depending on the dimension
of $Z$:
\item{{\sl (i)}} $\E_H(-1)=\O^{r+2}/\O(-1)^2$ and $\P(\E_H)$
has a contraction onto $\P^{r+1}$,
\item{{\sl (ii)}} $\E_H(-1)=(T\P^r(-1)\oplus\O(1))/\O$ and $\P(\E_H)$
admits a birational morphism onto a quadric $\Q^{2r-1}$, the variety
$\P(\E_H)$ is a blow-up of the quadric along a linear $\P^{r-1}$,
\item{{\sl (iii)}} $\E_H(-1)=\O^{r+1}/\O(-2)$
and $\P(\E_H)$ is contracted onto $\P^r$.

The remaining cases appearing in [PSW1, Thm.~1] are excluded:
decomposable bundles because of 4.2.2, the other ones because
they do not admit maps onto varieties of dimension emerging in cases
(i)---(iii) of 6.0.
\par
Note that above three cases are in one-to-one correspondence with
the cases (i)---(iii) from 6.0.
The variety $Z$ --- the target of the contraction $\f$
--- is therefore $\P^{r+1}$, $\Q^{2r}$ and $\P^r$, respectively.
If $\f$ is divisorial or a quadric bundle, we obtain a description
of $X$ (and therefore of $\E$) immediately --- see
6.3 and 6.1, respectively. If $\f$ is of type (i) then note that
$\E(-1)$ is spanned by $r+2$ sections, because $\O_{\P(\E(-1))}(1)=
\f^*\O(1)$, and therefore we have an
exact sequence
$$0\raa\H\raa\O^{r+2}\raa\E(-1)\raa0$$
with $\H$ a reflexive sheaf of rank 2. Since $\E(-1)$ restricted to
a hyperplane is $\O^{r+2}/\O(-1)^2$ it follows that $\H=\O(-1)^2$
and thus we have a description of $\E$ and of $X$.
\par
We summarize the result in the following

\proclaim Theorem 6.8.
Let $X$ be a Fano manifold of index $r$ and dimension $2r$.
Assume that $X$ is a projectivization of a sheaf $\E$, $p:X=\PE\ra Y$,
and assume moreover that $\E$ is not locally free.
Then one of the following holds (note that the top Chern class $c_{r+1}(\E)$
is equal to the number of singular fibres of $\E$):
\item{(i)} $Y\iso \P^{r+1}$, $X$ is an intersection of
two divisors of bidegree
$(1,1)$ in $\P^{r+1}\times\P^{r+1}$,
$$\E(-1)=\O^{r+2}/\O(-1)^2,\ \ c_{r+1}(\E)=r+2,$$
\item{(ii)} $Y\iso \P^{r+1}$ and $X$ is a blow-up of $\Q^{2r}$ along a linear
$\P^{r-1}\subset\Q^{2r}$,
$$\E(-1)=(T\P^{r+1}(-1)\oplus\O(1))/\O^2,\ \ c_{r+1}(\E)=2,$$
\item{(iii)} $Y\iso \P^{r+1}$, $X$ is a divisor of bidegree $(1,2)$ in
$\P^{r+1}\times\P^{r}$ and
$$\E(-1)=\O^{r+1}/\O(-2),\ \ c_{r+1}(\E)=2^{r+1},$$
\item{(iii')} $Y\iso \Q^{r+1}$, $X$ is a divisor of bidegree
$(1,1)$ in $\Q^{r+1}\times\P^r$ and
$$\E(-1)=\O^{r+1}/\O(-1),\ \ c_{r+1}(\E)=2.$$

\beginsection References.

\item{[ABW]} M.~Andreatta, E.~Ballico, J.~Wi\'sniewski, Vector bundles
and adjunction, Int.~J.~Math. {\bf 3} (1992), 331---340.
\item{[AW]} M.~Andreatta, J.~Wi\'sniewski, A note on non-vanishing
and applications, Duke Math. J. (to appear).
\item{[B]} C.~B\v anic\v a, Smooth reflexive sheaves,
Revue Romaine Math.~Pures Appl, {\bf 63} (1991), 571---573.
\item{[BC]} C.~B\v anic\v a, I.~Coand\v a,
Existance of rank 3 vector bundles on
homogeneous rational 3-folds, Manuscr. Math. {\bf 51} (1986), 121---143.
\item{[BSW]} M.~Beltrametti, A.~J.~Sommese and J.~Wi\'sniewski,
Results on varieties with many lines and their applications to
adjunction theory, in {\sl Complex Algebraic Varieties, Bayreuth 1990},
Lecture Notes in Math {\bf 1507}, Springer-Verlag 1992.
\item{[F1]} T.~Fujita, On polarized manifolds whose adjoint bundles
are not semipositive, in {\sl Algebraic Geometry, Sendai 1985},
Adv. Studies in Math. {\bf 10}, pp. 167--178, Kinokuniya 1987.
\item{[F2]} ------, On adjoint bundle of ample vector bundles,
in {\sl Complex Algebraic Varieties, Bayreuth 1990},
Lecture Notes in Math {\bf 1507}, Springer-Verlag 1992.
\item{[G]} A.~Grothendieck, EGA II, Publ.~Math.~IHES {\bf 8} (1961).
\item{[H1]} R.~Hartshorne, {\sl Algebraic Geometry},
Springer-Verlag 1977.
\item{[H2]} ------, Stable reflexive sheaves, Math.~Ann.~{\bf 254} (1980),
121---176.
\item{[H3]} ------, Ample vector bundles, Publ. Math. IHES
{\bf 29} (1966).
\item{[HH]} R.~Hartshorne, A.~Hirschowitz, Nouvelles curves de bon genre,
Math.~Ann.~{\bf 267} (1988), 353---367.
\item{[I]} P.~Ionescu, Generalized adjunction and applications,
Math.~Proc.~Camb.~Phil. Soc.~{\bf 99} (1988), 457---472.
\item{[KMM]} J.~Koll\'ar, Y.~Miyaoka, Sh.~Mori, Rational connectedness
and boudness of Fano manifolds, J.~Diff.~Geom.~{\bf 36} (1992), 765---779.
\item{[OSS]} Ch.~Okonek, M.~Schneider, H.~Spindler, {\sl Vector bundles
on complex projective spaces}, Progress in Math.~{\bf 3}, Birkh\"auser 1980.
\item{[PSW1]}  Th.~Peternell, M.~Szurek and J.~Wi\'sniewski,
Numerically effective vector bundles with small Chern classes,
in {\sl Complex Algebraic Varieties, Bayreuth 1990},
Lecture Notes in Math {\bf 1507}, Springer-Verlag 1992.
\item{[PSW2]} ----------,
Fano manifolds and vector bundles, Math.~Ann. (1992), 151---165.
\item{[S]} A.~J.~Sommese, On the adjunction
theoretic structure of projective varieties, in
{\sl Complex Analysis and Algebraic Geometry, G\"ottingen 1985},
Lecture notes in Math.~{\bf 1194}, 175---213, Springer-Verlag.
\item{[SW]} M.~Szurek, J.~Wi\'sniewski, On Fano manifolds which
are $\P^k$-bundles over $\P^2$, Nagoya Math.~J.~{\bf 120} (1990), 89---101.
\item{[W1]} J.~Wi\'sniewski, On contractions of Fano manifolds,
J.~reine Angew.~ Math.~{\bf 417} (1990), 141---157.
\item{[W2]} --------, Fano manifolds and quadric bundles,
to appear in Math.~Zeit.
\bigskip

\obeylines{
{\sl E.~Ballico: Dipartimento di Matematica, Universit\'a di Trento, 38050
Povo,
{\tt ballico@itnvax.science.unitn.it}
{\sl J.~Wi\'sniewski: Instytut Matematyki UW, Banacha 2, 02-097 Warszawa,
Poland
{\tt jarekw@plearn.edu.pl}}

\end